\documentclass[reqno,11pt]{article}

\usepackage{amsmath,latexsym,amssymb, amscd}

\def\sqr#1#2{{\vcenter{\hrule height.#2pt
        \hbox{\vrule width.#2pt height#1pt \kern#1pt
                \vrule width.#2pt}
        \hrule height.#2pt}}}
\def\square{\mathchoice\sqr64\sqr64\sqr{4}3\sqr{3}3}
\def\QED{\hfill$\square$\break}
\def\demo{\noindent{\bf Proof: }}

\DeclareMathOperator{\n}{\mathbf n}
\DeclareMathOperator{\q}{\mathbf q}

\DeclareMathOperator{\Spec}{\text{Spec }}

\DeclareMathOperator{\Ass}{\text{Ass }}

\newtheorem{Theorem}{\sc Theorem}[section]

\newtheorem{Proposition}[Theorem]{\sc Proposition}

\newtheorem{Remark}[Theorem]{\sc Remark}

\newtheorem{Question}[Theorem]{\sc Question}
\newtheorem{Construction}[Theorem]{\sc Construction of Example}
\newtheorem{Approximation}[Theorem]{\sc Approximation Technique}
\newtheorem{Comment}[Theorem]{\sc Comments}

\begin{document}

\baselineskip=13pt

\ \vspace{1.7in}

\noindent {\LARGE\bf Integral Closures of Ideals in Completions \\
\\
of Regular Local Domains}

\vspace{.25in}

\noindent WILLIAM \ HEINZER, \  Department of Mathematics, Purdue
University, West Lafayette, IN 47907-1395 USA  {\it E-mail}: {\tt
heinzer@math.purdue.edu}

\bigskip

\noindent CHRISTEL  \ ROTTHAUS, \  Department of Mathematics,
Michigan State University, East Lansing, MI 48824-1027 USA. {\it
E-mail}: {\tt rotthaus@math.msu.edu}

\bigskip

\noindent SYLVIA  \ WIEGAND, \  Department of Mathematics and
Statistics, University of Nebraska, Lincoln, NE 68588-0323 USA.
{\it E-mail}: {\tt swiegand@math.unl.edu}

\vspace{2.4cm}

\section{Abstract \hfill\break}

In this paper we exhibit  an example of a three-dimensional
regular local domain $(A, \n)$ having a height-two prime ideal $P$
with  the property that the extension $P\widehat A$  of $P$ to the
$\n$-adic completion $\widehat A$ of $A$ is not integrally closed.
We use a construction we have studied in earlier papers: For
$R=k[x,y,z]$, where $k$ is a field of characteristic zero and
$x,y,z$ are indeterminates over $k$, the example $A$ is an
intersection of the localization
 of the power series ring
$k[y,z][[x]]$ at the maximal ideal $(x,y,z)$ with the field
$k(x,y,z,f, g)$, where $f, g $ are elements of
$(x,y,z)k[y,z][[x]]$ that are algebraically independent over
$k(x,y,z)$. The elements $f, g$ are chosen in such a way that
using results from our earlier papers  $A$ is Noetherian and it is
possible to describe  $A$  as  a nested union of rings  associated
to $A$ that are  localized polynomial rings over $k$ in five
variables.

\bigskip

\section{Introduction and Background \hfill\break}

We are interested in the general question: What can happen 
in the completion of a `nice' Noetherian ring?   We are  examining 
this question as part of a project of 
constructing Noetherian and non-Noetherian integral domains using
power series rings. In this paper  as a continuation of that
project we display  an example of a three-dimensional regular
local domain $(A, \n)$  having a height-two prime ideal $P$ with
the property that the extension $P\widehat A$  of $P$ to the
$\n$-adic completion $\widehat A$ of $A$ is not integrally closed.
The ring $A$ in  the example is a nested union of
regular local domains of dimension five.

Let $I$ be an ideal of a commutative ring $R$ with identity. We
recall that an element $r \in R$ is {\it integral over} $I$ if
there exists a monic polynomial $f(x) \in R[x]$,
$f(x) = x^n + \sum_{i= 1}^na_ix^{n-i}$, where $a_i \in  I^i$ for
each $i, 1 \le i \le n$ and $f(r) = 0$. Thus $r \in R$ is
integral over $I$ if and only if 
$IJ^{n-1} = J^n$,   where $J = (I, r)R$ and $n$ is some positive 
integer. (Notice that $f(r) = 0$  implies 
$r^n = -\sum_{i=1}^na_ir^{n-i} \in IJ^{n-1}$ and this  
implies $J^n \subseteq IJ^{n-1}$.)
If $I \subseteq J$ are ideals and $IJ^{n-1} = J^n$, then 
$I$ is said to be a {\it reduction} of $J$.
The {\it integral closure} $\overline I$ of an ideal $I$ is the
set of elements of $R$ integral over $I$. If $I = \overline I$,
then $I$ is said to be {\it integrally closed}. 
It is well known that 
$\overline I$ is an integrally closed ideal. An ideal is 
integrally closed if and only if it is not a reduction of a 
properly bigger ideal. A prime ideal is always integrally
closed. An ideal is said to be {\it normal} if all the
powers of the ideal are integrally closed.

We were motivated to construct the  example given in
this paper by a question asked
by Craig Huneke as to whether there exists an analytically
unramified Noetherian local ring $(A, \n)$  having an integrally
closed ideal $I$ for which $I\widehat A$ is not integrally closed,
where $\widehat A$ is the $\n$-adic completion of $A$.
In Example \ref{const}, the ring $A$ is a 3-dimensional 
regular local domain and $I = P = (f,g)A$ is a prime
ideal of height two. Sam Huckaba asked if the ideal of
our example is a normal ideal. The answer is `yes'. Since
$f, g$ form a regular sequence and $A$ is Cohen-Macaulay, 
the powers $P^n$ of $P$ have no embedded associated 
primes and therefore are $P$-primary \cite[(16.F), p. 112]{M1},
\cite[Ex. 17.4, p.  139]{M2}. Since the powers of the maximal
ideal of a regular local domain are integrally closed, the 
powers of $P$ are integrally closed. Thus the Rees algebra
$A[Pt] = A[ft, gt]$ is a normal domain while the Rees
algebra $\widehat A[ft, gt]$ is not integrally closed.

A problem analogous to that considered here in the sense
that it also deals with the behavior of ideals under 
extension to completion is 
addressed by Loepp and Rotthaus in \cite{LR}. 
They construct nonexcellent local Noetherian domains 
to  demonstrate that tight closure 
need not commute with  completion.

\begin{Remark} {\rm
Without the assumption that $A$ is analytically unramified, there
exist examples even in dimension one where an integrally closed
ideal of $A$  fails  to extend to an integrally closed ideal in
$\widehat A$. If $A$ is reduced but analytically ramified, then
the zero ideal of $A$ is integrally closed, but its extension to
$\widehat A$ is not integrally closed. An example in  characteristic zero 
of a
one-dimensional Noetherian local domain that is analytically
ramified is given by  Akizuki in his  1935 paper  \cite{A}. 
An example in positive characteristic is given by F.K. Schmidt
\cite[pp.  445-447]{S}. 
Another  example due to Nagata is given in  \cite[Example 3,
pp.  205-207]{N}. (See also \cite[(32.2), p. 114]{N}.)}
\end{Remark}

\begin{Remark} \label{prel}
{\rm Let $R$ be a commutative ring and let $R'$ be an $R$-algebra.
We list  cases where extensions  to $R'$ of integrally closed
ideals of $R$ are again integrally closed. The $R$-algebra $R'$ is
said to be  {\it quasi-normal} if $R'$ is flat over $R$ and the
following condition $(N_{R,R'})$ holds: If $C$ is any $R$-algebra
and $D$ is a $C$-algebra in which $C$ is integrally closed, then
also $C \otimes_RR'$ is integrally closed in  $D\otimes_RR'$.}
\begin{enumerate}
\item
{\rm By \cite[Lemma 2.4]{L}, if $R'$ is an $R$-algebra satisfying
$(N_{R,R'})$ and $I$ is an integrally closed ideal of $R$, then
$IR'$ is integrally closed in $R'$.}

\item
{\rm Let $(A, \n)$ be a Noetherian local ring and let $\widehat A$
be the $\n$-adic completion of $A$. 
Since $A/\q \cong \widehat A/\q \widehat A$ for every $\n$-primary
ideal $\q$ of $A$, the $\n$-primary ideals of $A$ are in
one-to-one inclusion preserving correspondence with the 
$\widehat \n$-primary ideals of $\widehat A$. It follows that
an $\n$-primary ideal $I$ of $A$ is 
a reduction of a properly larger ideal of $A$  if and only if $I\widehat A$
is a reduction of a properly larger ideal of $\widehat A$.
Therefore an $\n$-primary ideal $I$ of $A$ is 
integrally closed if and 
only if $I\widehat A$ is integrally closed.}

\item
{\rm If $A$ is excellent,  then  the map $A \to
\widehat A$ is quasi-normal by \cite[(7.4.6) and (6.14.5)]{G}, and
in this case every integrally closed ideal of $A$ extends to an
integrally closed ideal of $\widehat A$.}

\item
{\rm If $(A, \n)$ is a local domain and $A^h$ is the
Henselization of $A$, then every integrally closed ideal of $A$
extends to an integrally closed ideal of $A^h$. This follows
because $A^h$ is a filtered direct limit of \'etale $A$-algebras
\cite[(iii), (i), (vii) and  (ix),  pp. 800- 801]{L}.}

\item
{\rm In general, integral closedness of ideals is a local
condition. 
Suppose $R'$ is an $R$-algebra  that is {\it locally normal}
in the sense that for every prime ideal $P'$ of $R'$, the local
ring $R'_{P'}$ is an integrally closed domain. Since principal
ideals of an integrally closed domain are integrally closed, the
extension to $R'$ of every principal ideal of $R$ is integrally
closed. In particular, if $(A, \n)$ is an analytically normal
Noetherian local domain, then every principal ideal of $A$ extends
to an integrally closed ideal of  $\widehat A$.}

\item
{\rm  If $R$ is an integrally closed domain, then for every
ideal $I$ and element $x$ of $R$ we have $\overline{xI} = x\overline I$.
If $(A, \n)$ is analytically normal and also a UFD, then
every height-one prime ideal of $A$ extends to an integrally
closed ideal of $\widehat A$. In particular if $A$ is a regular
local domain, then $P\widehat A$ is integrally closed for every
height-one prime $P$ of $A$. If $(A, \n)$ is a 2-dimensional
regular local domain, then every nonprincipal  integrally closed
ideal of $A$ has the form $xI$, where $I$ is an $\n$-primary
integrally closed ideal and $x \in A$. In view of item 2, every
integrally closed ideal of $A$ extends to an integrally closed
ideal of $\widehat A$ in the case where $A$ is a 2-dimensional
regular local domain. }

\item {\rm Suppose $R$ and $R'$ are Noetherian rings and assume
that $R'$ is a flat $R$-algebra. Let $I$ be an integrally closed
ideal of $R$. The flatness of $R'$ over $R$ implies every $P' \in
\Ass(R'/IR')$ contracts  in $R$ to some $P \in \Ass(R/I)$
\cite[Theorem 23.2]{M2}. Since a regular map is quasi-normal, if
the map $R \to R'_{P'}$ is regular for each $P' \in \Ass(R'/IR')$,
then $IR'$ is integrally closed.}
\end{enumerate}
\end{Remark}

\bigskip

\section{A non-integrally closed extension \hfill\break}

In the construction of the following example we make use of
results from \cite{HRW1}-\cite{HRW3}.

\bigskip

\begin{Construction} \label{const}
{\rm Let $k$ be a field  of characteristic zero and let $x, y$ and
$ z$ be indeterminates over $k$. Let $R :=  k[x, y, z]_{(x,y,z)}$
and let $R^*$ be the $(xR)$-adic completion of $R$. Thus $R^* =
k[y,z]_{(y, z )}[[x]]$, the formal power series ring in $x$ over
$k[y, z]_{(y, z)}$. }
\end{Construction}

Let $\alpha$ and $\beta$ be elements of $xk[[x]]$ which are
algebraically independent over $k(x)$. Set
$$
f = (y- \alpha)^2, \quad  g = (z- \beta)^2, \quad \text{ and } A =
k(x, y, z, f, g) \cap R^*.
$$
Then the $(xA)$-adic completion $A^*$ of $A$ is equal to $R^*$
\cite[Lemma 2.3.2, Prop. 2.4.4]{HRW2}.

In order to  better understand the structure of  $A$, we recall
some of the details of the construction of a nested union $B$  of
localized polynomial rings over $k$ in $5$ variables associated to
$A$.  (More details may be found in \cite{HRW3}.)

\begin{Approximation}
{\rm
 With $k, x, y, z,
f, g,  R$ and $R^*$ as in (3.1), Write
$$
f = y^2 + \sum^{\infty}_{j=1} b_{j}x^j, \qquad g = z^2 +
\sum^{\infty}_{j=1} c_{j}x^j,
$$
for some $b_j,\in k[y]$ and  $c_j \in k[z]$. There are natural
sequences $\{f_{r}\}_{r=1}^\infty$, $\{g_{r}\}_{r=1}^\infty$ of
elements in $R^*$, called the $r^{\text{th}}$ {\it endpiece}s for
$f$ and $g$ respectively which ``approximate"  $f$ and $g$. These
are defined for each $r \ge 1$ by:}
\end{Approximation}

$$\quad
f_{r}:=\sum_{j=r}^{\infty} (b_{j}x^j)/x^r ,\qquad
g_{r}:=\sum_{j=r}^{\infty} (c_{j}x^j)/x^r.
$$

For each $r \ge 1$, define $B_r$ to be $k[x, y, z, f_r, g_r]$
localized at the maximal ideal generated by $ (x, y, z, f_r - b_r,
g_r - c_r)$. Then  define $B = \bigcup_{r=1}^\infty B_r$. The
endpieces defined  here are slightly different from 
the notation used in \cite{HRW3}.  Also we are using here a
localized polynomial ring for  the base ring
$R$. With minor adjustments, however,     \cite[Theorem 2.2]{HRW3} applies
to our setup.

\begin{Theorem} \label{main}
Let  $A$ be the ring constructed in (3.1) and let $P = (f, g)A$,
where $f$ and $g$ are as in (3.1) and (3.2). Then
\begin{enumerate}
\item
$A = B$ is a three-dimensional regular local domain that is a
nested union of five-dimensional regular local domains.
\item
$P$ is a height-two prime ideal of $A$.
\item   If  $A^*$ denotes  the $(xA)$-adic
completion of $A$, then $A^* = k[y, z]_{(y, z)}[[x]]$ and $PA^*$
is not integrally closed.
\item
If  $\widehat A$ denotes  the completion of $A$ with respect to
the powers of the maximal ideal of $A$, then $\widehat A = k[[x,
y, z]]$ and $P\widehat A$ is not integrally closed.
\end{enumerate}
\end{Theorem}

\demo Notice that the polynomial ring
$k[x, y, z, \alpha, \beta] = k[x, y, z, y-\alpha, z-\beta]$ is
a free module of rank 4 over the polynomial subring
$k[x, y, z, f, g]$ since $f = (y-\alpha)^2$ and 
$g = (z-\beta)^2$. Hence the extension
$$
k[x, y, z, f, g] \to k[x, y, z, \alpha, \beta][1/x]
$$
is flat. Thus  item (1) follows from  \cite[Theorem 2.2]{HRW3}.

For item (2),  it suffices  to observe that $P$ has height two and
that,   for each positive integer $r$,    $P_r := (f, g)B_r$ is a
prime ideal of $B_r$. We have $f = xf_1 + y^2$ and $g =  xg_1 +
z^2$. It is clear that $(f, g)k[x, y, z, f, g]$ is  a height-two
prime ideal. Since $B_1$ is a localized polynomial ring over $k$
in the variables $x, y, z, f_1 -b_1, g_1 - c_1$, we see that
$$
P_1B_1[1/x] = (xf_1 + y^2, xg_1 + z^2)B_1[1/x]
$$
is a height-two prime ideal of $B_1[1/x]$. Indeed, setting $f = g
= 0$ is equivalent to setting $f_1 = -y^2/x$ and $g_1 = -z^2/x$.
Therefore the residue class ring $(B_1/P_1)[1/x]$ is isomorphic to
a localization of  the integral domain $k[x, y, z][1/x]$. Since
$B_1$ is Cohen-Macaulay and $f, g$ form a regular sequence, and
since $(x, f, g)B_1 = (x, y^2, z^2)B_1$ is an ideal of height
three, we see that $x$ is in no associated prime of $(f, g)B_1$
(see, for example \cite[Theorem 17.6]{M2}). Therefore $P_1 = (f,
g)B_1$ is a height-two prime ideal.

 For $r > 1$, there exist elements $u_r \in k[x, y]$ and $v_r
\in k[x, z]$ such that  $f = x^rf_r  + u_rx + y^2$ and $g = x^rg_r
+ v_rx + z^2$. An argument similar to that given above  shows that
$P_r = (f, g)B_r$ is a height-two prime of $B_r$. Therefore $(f,
g)B$ is a height-two prime of $B=A$.

For items 3 and 4, $R^*=B^*=A^*$ by Construction \ref{const} and
it follows that $\widehat A = k[[x, y, z]]$. To see that $PA^* =
(f,g)A^*$ and $P\widehat A = (f, g)\widehat A$ are not integrally
closed, observe that $\xi := (y-\alpha)(z-\beta)$ is integral over
$PA^*$ and $P\widehat A$  since $\xi^2 = fg \in P^2$. On the other
hand, $y-\alpha$ and $z-\beta$ are nonassociate prime elements in
the local unique factorization domains  $A^*$ and $\widehat A$. An
easy computation shows that $\xi \not\in P\widehat A$.  Since
$PA^* \subseteq P\widehat A$, this completes the proof. \QED

\smallskip

\begin{Remark}
{\rm In a similar manner it is possible to construct for each
integer $d \ge 3$ an example of a $d$-dimensional regular local
domain $(A, \n)$ having a prime ideal $P$ of height $h := d-1$
such that $P\widehat A$ is not integrally closed. Indeed, let $k$
be a field of characteristic zero and let $x, y_1, \ldots, y_h$ be
indeterminates over $k$. Let $\alpha_1, \ldots, \alpha_h \in
xk[[x]]$ be algebraically independent over $k(x)$. For each $i$
with $1 \le i \le h$, define $f_i = (y_i - \alpha_i)^h$.
Proceeding in a manner similar to what is done in (3.1) we obtain
a $d$-dimensional regular local domain $A$ and a prime ideal $P =
(f_1, \ldots, f_h)A$ of height $h$ such that the $y_i- \alpha_i
\in \widehat A$. Let $\xi = \prod_{i=1}^h(y_i - \alpha_i)$. Then
$\xi^h = f_1 \cdots f_h  \in P^h$ implies $\xi$ is integral over
$P\widehat A$, but using that $y_1-\alpha_1, \ldots, y_h
-\alpha_h$ is a regular sequence in $\widehat A$, we see that $\xi
\not\in P\widehat A$.}
\end{Remark}

\smallskip

\section{Comments and Questions  \hfill\break}

In connection with Theorem \ref{main} it is natural to ask the
following question.

\begin{Question} \label{4.1}
{\rm For $P$ and $A$ as in Theorem \ref{main}, is $P$ the only
prime of $A$ that does not extend to an integrally closed ideal of
$\widehat A$? }
\end{Question}

\begin{Comment} \label{Com}
{\rm In relation to the  example given in Theorem \ref{main}
 and to Question \ref{4.1}, we have the
following commutative diagram, where all the maps shown are the
natural inclusions:}

\begin{equation}
\CD B=A @>{\gamma_1}>> A':=k(x,y,z,\alpha,\beta)\cap R^*
@>{\gamma_2}>> R^*=A^*
\\
 @A{\delta_1}AA        @A{\delta_2}AA        @A{\psi}AA           \\
S:=k[x,y,z,f,g] @>{\varphi}>> T := k[x,y,z,\alpha,\beta]@= T \\
\endCD
\end{equation}

{\rm Let $\gamma=\gamma_2\cdot \gamma_1$. Referring to the diagram
above, we observe the following:}
\begin{enumerate}
\item{\rm  The discussion in \cite[bottom  p. 668 to top p. 669]{HRW2} 
implies that \cite[Thm. 3.2]{HRW2} applies to  the setting of Theorem~3.3.
By \cite[Prop. 4.1 and Thm. 3.2]{HRW2}, $A'[1/x]$ is
a localization of $T$. By Theorem 3.3 and \cite[Thm 3.2]{HRW2}, 
$A[1/x]$ is a localization of $S$.
Furthermore, by  \cite[Prop. 4.1]{HRW2}
$A'$ is excellent. (Notice, however, that $A$ is not excellent
since there exists  a prime ideal $P$ of $A$ such that  $P\widehat
A$ is not integrally closed.) The excellence of $A'$ implies that
if $Q^*\in\Spec A^*$ and $x\notin Q^*$, then the map $\psi_{Q^*}:
T\to A^*_{Q^*}$ is regular \cite[(7.8.3 v)]{G}.
\item
Let $Q^* \in \Spec A^*$ be such that $x\notin Q^*$ and let $\q' =
Q^* \cap T$. By \cite[Theorem 32.1]{M2} and Item 1 above, if
$\varphi_{\q'}: S \to T_{\q'}$ is regular, then $\gamma_{Q^*}:
A\to A^*_{Q^*}$ is regular.
\item
Let $I$ be an ideal of $A$. Since $A'$ and  $A^*$ are excellent
and both have completion $\widehat A$, Remark \ref{prel}.3 shows
that  the ideals $IA'$, $IA^*$ and $I\widehat A$ are either all
integrally closed or all fail to be integrally closed.
\item
The Jacobian ideal of the extension $\varphi : S = k[x, y, z, f,
g] \to T = k[x, y, z, \alpha, \beta]$ is the ideal of $T$
generated by the determinant of the matrix
$$
{ \cal J}:= \left(\begin{array}{cc} \frac{\partial f}{\partial \alpha} &
\frac{\partial g}{\partial \alpha}
\\[.2cm]
\frac{\partial f}{\partial \beta} & \frac{\partial g}{\partial
\beta}
\end{array}\right).
$$
Since the characteristic of the field $k$ is zero, this ideal is
$(y-\alpha)(z-\beta)T$.}
\end{enumerate}
\end{Comment}

In Proposition \ref{4.3}, we relate the behavior of integrally
closed ideals in the extension $\varphi :S \to T$ to the behavior
of integrally closed ideals in the extension $\gamma : A \to A^*$.

\begin{Proposition} \label{4.3}
With the setting of Theorem \ref{main} and Comment \ref{Com}.2,
let $I$ be an integrally closed ideal of $A$ such that $x \not\in
Q$ for each $Q \in \Ass(A/I)$. Let $J  = I \cap S$. If $JT$ is
integrally closed (resp. a radical ideal) then $IA^*$ is
integrally closed (resp.  a radical ideal).
\end{Proposition}
\demo Since the map $A \to A^*$ is flat, $x$ is not in any
associated prime of $IA^*$. Therefore $IA^*$ is contracted from
$A^*[1/x]$ and it suffices to show $IA^*[1/x]$ is integrally
closed (resp. a radical ideal). 
Our hypothesis implies $I = IA[1/x] \cap A$. 
By Comment \ref{Com}.1,  $A[1/x]$ is a localization of $S$. 
Thus every ideal of $A[1/x]$ is the extension of its 
contraction to $S$. It follows that $IA[1/x] = JA[1/x]$. 
Thus $IA^*[1/x] = JA^*[1/x]$.

Also by  Comment \ref{Com}.1,  the map $T \to A^*[1/x]$ is regular.
If  $JT$ is integrally closed, then  Remark
\ref{prel}.7 implies that $JA^*[1/x]$ is integrally closed.
If $JT$ is a radical ideal, then  
the  regularity of the map $T \to A^*[1/x]$ implies the $JA^*[1/x]$ is a
radical ideal.  \QED

\begin{Proposition}
With the setting of Theorem \ref{main} and Comment \ref{Com}, let
$Q\in \Spec A$ be such that $Q\widehat A$ (or equivalently $QA^*$)
is not integrally closed. Then
\begin{enumerate}
\item $Q$ has height two and $x \not\in Q$.
\item There exists a minimal prime $Q^*$ of $QA^*$
such that with $\q' = Q^* \cap T$, the map $\varphi_{\q'}: S \to
T_{\q'}$ is not regular.
\item $Q$ contains $f = (y-\alpha)^2$ or $g = (z-\beta)^2$.
\item
$Q$ contains no element that is a regular parameter of $A$.
\end{enumerate}
\end{Proposition}
\demo By Remark \ref{prel}.6, the height of $Q$ is two. Since
$A^*/xA^* = A/xA = R/xR$, we see that $x \not\in Q$. This proves
item 1.

By Remark \ref{prel}.7, there exists a minimal prime $Q^*$ of
$QA^*$ such that $\gamma_{Q^*}: A\to A^*_{Q^*}$ is not regular.
Thus item 2 follows from Comment \ref{Com}.2.

For  item 3,  let $Q^*$ and $\q'$ be as in item 2. 
Since $\gamma_{Q^*}$ is not regular it is not essentially
smooth \cite[6.8.1]{G}. 
By \cite[ (2.7)]{HRW3},
$(y-\alpha)(z-\beta) \in \q'$. Hence $f = (y - \alpha)^2$ or $g =
(z-\beta)^2$ is in $\q'$ and thus in $Q$. This proves item 3.

Suppose $w \in Q $ is a regular parameter for $A$. Then $A/wA$ and
$A^*/wA^*$ are two-dimensional regular local domains.  By Remark
\ref{prel}.6, $QA^*/wA^*$ is integrally closed, but this implies
that $QA^*$ is integrally closed, which contradicts our
 hypothesis that $QA^*$ is
not integrally closed. This proves item~4. \QED

\begin{Question}
{\rm In the setting of Theorem \ref{main} and Comment \ref{Com},
let $Q \in \Spec A$ with $x \notin Q$ and let $\q = Q\cap S$. If
$QA^*$ is integrally closed, does it follow that $\q T$ is
integrally closed?}
\end{Question}

\begin{Question}
{\rm In the setting of Theorem \ref{main} and Comment \ref{Com},
if  a prime ideal $Q$ of $A$ contains $f$ or
$g$, but not both,  and does not contain a regular parameter
of $A$, does it follow that  $QA^*$  is integrally closed ?}
\end{Question}

In Example \ref{const}, the three-dimensional regular local domain
$A$ contains  height-one prime ideals $P$ such that $\widehat
A/P\widehat A$ is not reduced. This motivates us to ask:

\begin{Question}
{\rm Let $(A, \n)$ be a three-dimensional regular local domain and
let $\widehat A$ denote the $\n$-adic completion of $A$.  If for
each height-one prime $P$ of $A$, the extension $P\widehat A$ is a
radical ideal, i.e., the ring $\widehat A/P\widehat A$ is reduced,
does it follow that $P\widehat A$ is integrally closed for each $P
\in \Spec A$?}
\end{Question}

\begin{center}
{\Large\bf Acknowledgments}
\end{center}

  The authors are grateful for
  the hospitality and cooperation of Michigan State, Nebraska
   and Purdue, where several work sessions on this research were
   conducted.  The authors also thank the referee for 
helpful suggestions and comments about the paper.

\end{document}